\magnification=1000
\hsize=11.7cm
\vsize=18.9cm
\lineskip2pt \lineskiplimit2pt
\nopagenumbers

\hoffset=-1truein
\voffset=-1truein

\advance\voffset by 4truecm
\advance\hoffset by 4.5truecm

\newif\ifentete

\headline{\ifentete\ifodd   \count0
      \rlap{\head}\hfill\tenrm\llap{\the\count0}\relax
    \else
        \tenrm\rlap{\the\count0}\hfill\llap{\head} \relax
    \fi\else
\global\entetetrue\fi}

\def\entete#1{\entetefalse\gdef\head{#1}}
\entete{}

\input amssym.def
\input amssym.tex

\def\-{\hbox{-}}
\def\.{{\cdot}}
\def\O{{\cal O}}

\def\F{{\cal F}}

\def\int{\frak i\frak n\frak t}

\def\qq{\quad{\rm and}\quad}

\def\too{\longrightarrow}

\font\large=cmr10  scaled \magstep 2
 2
 2
 2
\font\cds=cmr7
\font\cdt=cmti7

\centerline{\large On blocks with trivial source simple modules}

\bigskip
\centerline{\bf Lluis Puig }

\medskip
\centerline{\cds   CNRS, Institut de Math\'ematiques de Jussieu}
\smallskip
\centerline{\cds   6 Av
Bizet, 94340 Joinville-le-Pont, France}
\smallskip
\centerline{\cds E-mail: puig@math.jussieu.fr}

\bigskip
\centerline{\bf Yuanyang Zhou}

\medskip
\centerline{\cds Department of Mathematics and Statistics}
\smallskip
\centerline{\cds Central China Normal University, Wuhan, 430079, P.R. China}
\smallskip
\centerline{\cds E-mail: zhouyy74@163.com}

\vskip 0.5cm
\noindent
{\bf £1. Introduction}
\bigskip
£1.1. In [3] Danz and K\"ulshammer, investigating the simple modules for the large Mathieu groups, have found two blocks with noncyclic defect groups of order 9 where all the simple modules have trivial sources
and whose source algebras are isomorphic to the source
algebras of the corresponding blocks of their
 {\it inertial subgroups\/} [3, Theorems~4.3 and~4.4]{\footnote{\dag}{\cds As a matter of fact, from [12,~Corollary~3.6]  one easily may find infinitely many examples of such blocks.}}.

 \medskip
 £1.2. In their Introduction they note that, in general, any simple module with a trivial source determines an
 Alperin's {\it weight\/} [1] --- for instance, this follows   from [8, Proposition~1.6] --- and therefore, in a block with Abelian defect groups and all the simple modules with trivial  sources, Alperin's conjecture in [1] forces a canonical bijection between the sets of isomorphism classes of simple  modules of the block and of the corresponding block
 of its {\it inertial subgroup\/}. From this remark, they raise the question whether, behind this bijection, it should be a
true  Morita equivalence between both blocks.

\medskip
£1.3. Recently, Zhou proved that, in a suitable inductive context, the answer is in the affirmative [18, Theorem B]; our purpose here is to prove the same fact without any hypothesis on the defect group. In order to explicit
our result we need some notation;  let $p$ be a prime number, $k$ an algebraically closed field of
characteristic~$p\,,$ $G$ a finite group, $b$ a primitive idempotent of the center $Z(kG)$ of the group
algebra of $G$ --- for short, a {\it block\/} of $G$ --- and $P_\gamma$ a {\it defect pointed group\/} of
$b\,;$ that is to say, $P$ is a {\it defect group\/} of this block in Brauer's terms and $\gamma$ is a conjugacy class
 of primitive idempotents~$i$ in $(kGb)^P$ such that ${\rm Br}_P(i)\not= 0\,;$ here, ${\rm Br}_P$ denotes the
 usual {\it Brauer homomorphism\/}
  $${\rm Br}_P : (kG)^P\too (kG)(P) = (kG)^P\Big/\sum_Q (kG)^P_Q\cong kC_G (P)
  \eqno £1.3.1\phantom{.}$$
  where $Q$ runs over the set of proper subgroups of $P\,.$ Recall that the {\it $P\-$interior algebra\/}
  $(kG)_\gamma = i(kG)i$ is called a {\it source algebra\/} of $b$ and that its underlying $k\-$algebra is {\it Morita equivalent\/} to $kGb$ [8, Definition~3.2 and Corollary~3.5].

  \medskip
  £1.4. If~$G'$ is a second finite group and $b'$ a block of  $G'$ admitting the {\it same\/} defect group $P\,,$ it follows from [13, Corollary~7.4 and Remark~7.5]
  that the source algebras of $b$ and $b'$ are isomorphic --- as {\it $P\-$interior algebras\/} --- if and only if
  the categories of finitely generated $kGb\-$ and $kG'b'\-$modules are equivalent to each other {\it via\/} a
  $kGb\otimes_k kG'b'\-$module admitting a $P\times P\-$stable basis, a fact firstly proved by Leonard Scott [17,~Lemma]{\footnote{\dag}{\cds Strictly speaking, in
  [17, Lemma] Scott only considers the case where the {\cdt block algebras}
 {\cdt kGb} and {\cdt kG'b'} are isomorphic.}}; in this case, we simply say that the blocks $b$ and $b'$ are
 {\it identical\/}. More generally, we say that  the blocks $b$ and $b'$ are  {\it stably identical\/}  if
  the categories of finitely generated $kGb\-$ and $kG'b'\-$modules are {\it stably equivalent\/} to each other
  --- namely, equivalent to each other up to projective modules --- {\it throughout\/} a  $kGb\otimes_k kG'b'\-$module admitting a $P\times P\-$stable basis.

 \medskip
 £1.5. Set $N = N_G (P_\gamma)$ ---  often called the {\it inertial subgroup\/} of $b$ --- and denote by $e$ the block of $C_G (P)$ determined by  {\it the local point\/} $\gamma$ (cf.~£1.3.1). Recall that $e$ is also a block of $N$ and that
$k\bar C_G(P)\bar e$ is a simple $k\-$algebra, where we set $\bar C_G (P) = C_G (P)/Z(P)$ and denote by
$\bar e$ the image of~$e$ in~$k\bar C_G (P)\,;$ then, the action of $N$ on the simple $k\-$algebra
$k\bar C_G(P)\bar e$ determines a central $k^*\-$extension $\hat E$ of $E = N/P\.C_G (P)$ ---
often called the {\it inertial quotient\/} of~$b\,.$
Setting $\hat L = P\rtimes \hat E^\circ\,$ for a lifting of the canonical homomorphism $\hat E\to {\rm Out}(P)$ to ${\rm Aut}(P)\,,$ it follows from [11, Proposition~14.6] that the corresponding {\it twisted\/} group algebra  $k_*\hat L$ is isomorphic to a source algebra of~the~block~$e$ of~$N\,.$

\medskip
£1.6.  Recall that a {\it Brauer $(b,G)\-$pair} $(Q,f)$ is formed by a $p\-$subgroup $Q$ of $G$ such that
${\rm Br}_Q(b) \not= 0$ and by a block $f$ of $C_G(Q)$ fulfilling ${\rm Br}_Q(b)f = f$ [2,~Definition~1.6];
note that $f$ is also a block for any subgroup $H$ of $N_G(Q,f)$ containing $C_G(Q)\,.$ Thus,
$(P,e)$ is a  Brauer $(b,G)\-$pair and, as a matter of fact, there is $x\in G$ such that [2,~Theorem~1.14]
$$(Q,f)\i (P,e)^x
\eqno £1.6.1.$$
Then, the {\it Frobenius category\/}   $\F_{\!(b,G)}$ of $b$ [16,~3.1] is the category where the objects are the
 Brauer $(b,G)\-$pairs $(Q,f)$ and the morphisms are the homomorphisms between the corresponding $p\-$groups
 induced by the {\it inclusion\/} between Brauer $(b,G)\-$pairs and the $G\-$conjugation.

 \medskip
 £1.7. For short,  let us say that the block $b$ is {\it inertially controlled\/} whenever  the  Frobenius categories $\F_{\!(b,G)}$ and $\F_{\!\hat L}$ are equivalent to each other --- note that the unity element is the unique block of $\hat L$ and  we omit to mention it; moreover, since $k_*\hat L$  is isomorphic to a source algebra of the block $e$ of~$N\,,$  the  Frobenius
 categories $\F_{\!(e,N)}$ and $\F_{\!\hat L}$ are always equivalent to each other,  so that $e$ is always   {\it inertially controlled\/}.  Similarly, let us say that $b$ is a  {\it block of $G$ with trivial simple modules\/} if all the  simple $kGb\-$modules have trivial sources.

\bigskip
\noindent
{\bf Theorem~£1.8.} {\it With the notation above, the source algebra  $(kG)_\gamma$ of the block $b$ of $G$ is isomorphic to $k_*\hat L$ if and only if  the block $b$ of $G$ is inertially controlled and, for any Brauer $(b,G)\-$pair $(Q,f)$ contained in $(P,e)\,,$  $f$ is a block of $C_G (Q)\.N_P (Q)$ with  trivial source simple modules.\/}

\bigskip
£1.9. The main tools in proving this result are the   Linckelmann's Equi-valence Criterion on {\it stable equivalences\/} [7, Proposition~2.5], the {\it strict semicovering\/} homomorphisms that we recall in \S3 below, and the general criterion
 on {\it stable equivalences\/} in~[13,~Theorem~6.9], which in our context is summarized by the following result.

\bigskip
\noindent
{\bf Theorem~£1.10.} {\it With the notation above, the blocks $b$ of $G$ and $e$ of $N$ are stably identical if and only if, for any nontrivial Brauer $(b,G)\-$pair  $(Q,f)$ contained in $(P,e)\,,$ the block  $f$ of $C_G (Q)$ admits $C_P(Q)$ as a  defect $p\-$subgroup and a source algebra isomorphic to $k_*\big(C_{\hat L}(Q)\big)\,.$\/}

\bigskip
£1.11. Note that $C_{\hat E}(Q)$ acts faithfully on $C_P(Q)$ since any ($p'\-$)subgroup  of~$C_{\hat E}(Q)$ acting trivially on $C_P(Q)$ still acts trivially on $P$ [5, Ch.~5, Theorem~3.4], and that  we actually have
 $$C_{\hat L}(Q)\cong C_P(Q)\rtimes C_{\hat E^\circ}(Q)
 \eqno£1.11.1.$$
Moreover,  if the defect group $P$ is Abelian then, for any  Brauer $(b,G)\-$pair $(Q,f)$  contained in $(P,e)\,,$ $P$
is clearly a defect group of the block $f$ of $C_G(Q)\,.$ Finally, although we only work over $k\,,$ Lemma~7.8 in [10] allows us to lift all the isomorphisms between {\it block source algebras\/} over $k$ above to the corresponding block source algebras over a complete discrete valuation ring $\O$ of characteristic zero having the {\it residue field\/}
 $k\,.$

\bigskip
\bigskip
\noindent
{\bf £2. Notation and quoted results}
\bigskip

£2.1. Let $A$ be a finitely dimensional $k\-$algebra; we denote by
$1_A$ the unity element of $A$ and by $A^*$ the multiplicative group
of $A$. An algebra homomorphism $f$ from $A$ to another finitely
dimensional $k\-$algebra $A'$ is not necessarily unitary and we say
that $f$ is an {\it embedding\/} whenever
$${\rm Ker}(f) = \{0\}\qq {\rm Im}(f) = f(1_A)A'f(1_A)
\eqno £2.1.1.$$  Following Green, a {\it $G\-$algebra\/} is a finitely  dimensional $k\-$algebra $A$  endowed with a $G\-$action;  recall that, for any subgroup $H$ of $G\,,$ a {\it point\/} $\alpha$ of $H$ on $A$ is  an $(A^H)^*\-$conjugacy class of primitive idempotents of $A^H$ and the pair $H_\alpha$ is called a  {\it pointed group\/} on $A$ [8, 1.1]; we denote
 by~$A(H_\alpha)$ the {\it simple quotient\/} of $A^H$ determined
 by~$\alpha\,.$ A second pointed group $K_\beta$ on $A$ is {\it contained\/}
in~$H_\alpha$ if $K\i H$ and, for any $i\in\alpha\,,$ there is $j\in \beta$ such that~[8,~1.1]
$$ij = j = ji
\eqno £2.1.2.$$

\medskip
£2.2. Following Brou\'e, for any $p\-$subgroup $P$ of $G$ we consider the {\it Brauer quotient\/} and the {\it Brauer homomorphism\/}
$${\rm Br}^A_P : A^P\too A (P) =  A^P\Big/\sum_Q A^P_Q
\eqno £2.2.1,$$
where $Q$ runs over the set of proper subgroups of $P$ and $A_Q^P$ is the ideal formed by the sums $\sum_u a^u$ where $a$ runs over $A^Q$ and $u\in P$ over a set of representatives for $P/Q\,;$ we call  {\it local\/} any point~$\gamma$ of $P$
on $A$ not contained in   ${\rm Ker(Br}^A_P)$ [8, 1.1]. Let us say that $A$ is a {\it $p\-$permutation $G\-$algebra\/} if a Sylow $p\-$subgroup of $G$ stabilizes a basis of $A\,;$ in this case, recall that if $P$
 is a $p\-$subgroup of $G$ and $Q$ a normal subgroup of $P$
then the corresponding Brauer homomorphisms induce a $k\-$algebra isomorphism [2, Proposition~1.5]
$$\big(A(Q)\big)(P/Q)\cong A(Q)
\eqno £2.2.2.$$
 Obviously, the group algebra $A=kG$ is a
$p$-permutation $G$-algebra and the composition of the inclusion $k
C_G(Q)\subset A^Q$ with ${\rm Br}^A_Q$ is an isomorphism
which allows us to identify $k C_G(Q)$ with $A(Q)\,;$ then any local point
$\delta$ of $Q$ on $k G$ determines a block $b_\delta$ of $k C_G(Q)$
such that $b_\delta{\rm Br}_Q^{k G}(\delta)={\rm Br}_Q^{k
G}(\delta)$.

\medskip
£2.3. We are specially interested in the $G\-$algebras $A$ endowed with a group
homomorphism $\rho\,\colon G\to A^*$ inducing the action of $G$ on $A$ --- called {\it $G\-$ interior algebras\/}.
 In this case, for any pointed group $H_\alpha$ on $A$ and any $i\in \alpha\,,$ the subalgebra $A_\alpha = iAi$ has a  structure of  {\it $H\-$interior algebra\/} mapping~$y\in H$ on $\rho (y)i = i\rho(y)\,;$ moreover, setting $x\.a\.y = \rho(x)a\rho(y)$
for any $a\in A$ and any $x,y\in G\,,$ a  $G\-$interior algebra homomorphism
from~$A$ to another   $G\-$interior algebra~$A'$ is a $G\-$algebra homomorphism $f\,\colon A\to A'$ fulfilling
$$f(x\.a\.y) = x\.f(a)\.y
\eqno £2.3.1.$$
We also consider the {\it mixed\/} situation of an {\it $H\-$interior $G\-$algebra
$B$\/} where $H$ is a subgroup of $G$ and $B$ is a $G\-$algebra endowed with a {\it compatible\/} $H\-$interior algebra structure, in such a way that the
$kG\-$module $B\otimes_{kH} kG$ endowed with the product
$$(a\otimes x).(b\otimes y) = ab^{x^{-1}}\otimes xy
\eqno £2.3.2,$$
for any $a,b\in B$ and any $x,y\in G\,,$ and with the group homomorphism mapping $x\in G$ on $1_B\otimes x$  becomes a $G\-$interior algebra --- simply noted $B\otimes_H G\,.$ For instance, for any $p\-$subgroup $P$ of $G\,,$
$A(P)$ is a $C_G (P)\-$interior $N_G (P)\-$algebra.
\eject

\medskip
£2.4. In particular, if $H_\alpha$ and $K_\beta$ are two pointed groups on $A\,,$ we say that an
injective group homomorphism $\varphi\,\colon K\to H$ is an
{\it $A\-$fusion from $ K_\beta$ to~$H_\alpha$\/} whenever there is
a  $K\-$interior algebra~{\it embedding\/}
$$f_{\varphi} : A_\beta\too {\rm Res}^{H}_{K} (A_\alpha)
\eqno £2.4.1\phantom{.}$$
such that the inclusion $A_\beta\i A$ and the composition of $f_{\varphi}$
with the inclusion $A_\alpha\i A$ are $A^*\-$conjugate; we denote by $F_A
( K_\beta,H_\alpha)$ the set of $H\-$conjugacy classes of  $A\-$fusions from $ K_\beta$ to~$H_\alpha$
and we write $F_A (H_\alpha)$ instead of $F_A(H_\alpha,H_\alpha)\,.$ If $A_\alpha = iAi$ for $i\in \alpha\,,$
 it follows from [9, Corollary~2.13] that we have a group homomorphism
$$F_A (H_\alpha)\too N_{A_\alpha^{^*}}(H\.i)\big/H\.(A_\alpha^H)^*
\eqno £2.4.2.$$

\medskip
£2.5. Let $b$ be a block of $G\,;$ then $\alpha = \{b\}$ is a {\it
point\/} of $G$ on $kG$ and we let $P_\gamma$ be a local pointed
group contained in $G_\alpha$ which is maximal with respect to the
inclusion of pointed groups; namely $P_\gamma$ is a {\it defect
pointed group} of $b$. Note that, for any $p\-$subgroup $Q$ of~$G$
and any subgroup $H$ of $N_G (Q)$ containing $Q\,,$ we have
$${\rm Br}_Q\big((kG)^H\big) = (k C_G(Q))^H
\eqno £2.5.1;$$
thus, we have an injection from the set of points of $H$ on $kC_G (Q)$ to the set of points of $H$ on $kG$ such that the corresponding points $\beta^\circ$ and $\beta$ fulfill ${\rm Br}_Q^{k G} (\beta)={\rm Br}_Q^{k C_G(Q)}
(\beta^\circ)\,;$ moreover, this injection preserves the localness and the inclusion of pointed groups [16, 1.19].
In particular, if $P$ is Abelian and $Q_\delta$ is a local pointed group on $kG$ contained
in~$P_\gamma\,,$ a point $\mu$ of $C_G (Q)$ on $kG$ fulfilling
$$Q_\delta\i P_\gamma\i C_G(Q)_\mu
\eqno~£2.5.2\phantom{.}$$
 is the {\it unique\/} point determined by the block $b_\delta$ of $C_G (Q)$ and therefore $P$ is a defect group of this block
  (cf.~£1.8).

\medskip
£2.6. Set $e=b_\gamma$ and $N=N_G(P_\gamma)\,;$ thus, $e$ is a block
of $N\,,$ it determines a point $\nu$ of $N$ on $kG$ (cf.~£2.5) and $P$ is a defect group of this block; moreover, we have (cf.~£1.3.1)
$$(kN)(P)\cong kC_N (P) = kC_G(P)\cong (kG)(P)
\eqno £2.6.1,$$
there is a local point $\hat\gamma$ of $P$ on $kN\i kG$ such that
${\rm Br}_P(\hat \gamma)  = {\rm Br}_P(\gamma)$ and it follows from [4, Proposition~4.10] that, for any $\hat\imath\in \hat\gamma$ and any
$\ell\in \nu\,,$ the idempotent $\hat\imath\ell$ belongs to $\gamma$ and that the multiplication by $\ell$ defines a unitary $P\-$interior algebra homomorphism (cf.~£1.5)
$$k_*\hat L\cong (kN)_{\hat\gamma}\too (kG)_{\gamma}
\eqno £2.6.2\phantom{.}$$
which is actually a {\it direct injection\/} of $k(P\times P)\-$modules.
\eject

\medskip
£2.7. For any pair of {\it local pointed groups\/} $Q_\delta$ and $R_\varepsilon$ on $kG\,,$ we denote
by~$E_G (R_\varepsilon,Q_\delta)$ the set of $Q\-$conjugacy classes of group homomorphisms
$\varphi\,\colon R\to Q$ induced the conjugation by some $x\in G$ fulfilling $R_\varepsilon\i (Q_\delta)^x\,,$
and write $E_G(Q_\delta)$ instead of $E_G(Q_\delta,Q_\delta)\,;$
it follows from [9, Theorem~3.1] that
$$E_G (R_\varepsilon,Q_\delta) = F_{kG}(R_\varepsilon,Q_\delta)
\eqno £2.7.1\phantom{.}$$
and if $P_\gamma$ contains $Q_\delta$ and $R_\varepsilon$ then they can be considered as local pointed groups
on $(kG)_\gamma$ and it follows from [9, Proposition~2.14] that
$$E_G (R_\varepsilon,Q_\delta) = F_{kG}(R_\varepsilon,Q_\delta) = F_{(kG)_\gamma}(R_\varepsilon,Q_\delta)
\eqno £2.7.2.$$
In particular, it is clear that $N_G(Q_\gamma)/Q\.C_G (Q) \cong E_G(Q_\delta)$ and the action  of~$N_G (Q_\delta)$ on the simple $k\-$algebra $(kG)(Q_\delta)$ (cf.~£2.1) determines a central $k^*\-$extension $\hat E_G(Q_\delta)$ of
$E_G(Q_\delta)\,.$

\medskip
£2.8. Recall that a Brauer $(b,G)\-$pair $(Q,f)$ is called {\it selfcentralizing\/} if, setting $\bar C_G (Q) =C_G (Q)/Z(Q)$
and denoting by $\bar f$ the image of $f$ in $k\bar C_G (Q_\delta)\,,$ the  $k\-$algebra  $k\bar C_G (Q)\bar f$ is simple  [14,~1.6], so that  $k\bar C_G (Q)\bar f\cong (kG)(Q_\delta)$ for a local point $\delta$ of $Q$ on $kG$ clearly determined
by $f\,;$ we also say that $Q_\delta$ is a {\it selfcentralizing pointed group\/} on $kG\,;$ thus we have a bijection,
which preserves {\it  inclusion\/} and  {\it $G\-$conjugacy\/}, between the sets of  selfcentralizing pointed groups on $kGb$ and of selfcentralizing Brauer $(b,G)\-$pairs. Moreover, according to [14, Theorem~A.9],  an {\it essential pointed group\/} on $kG$ is a selfcentralizing pointed group $Q_\delta$ on $kG$ fulfilling the following condition
\smallskip
\noindent
£2.8.1\quad {\it $E_G (Q_\delta)$ admits a proper subgroup $M$ such that $p$ divides $\vert M\vert$
and does not divide $\vert M\cap M^\sigma\vert$ for any $\sigma\in E_G (Q_\delta) -M\,.$\/}
\smallskip
\noindent
Then,  from [14,~Corollary~A.12] and [16,~Corollary~5.14], it is not difficult to prove
that the block $b$ of $G$ is {\it inertially controlled\/} (cf.~£1.7) if and only if there are {\it no essential pointed groups\/} on $kGb\,;$ thus, if the defect group $P$ is Abelian  the block~$b$ of $G$ is inertially controlled.

\bigskip
\noindent
{\bf Lemma~£2.9.} {\it With the notation above, the block $b$ of $G$ is inertially controlled  if and only if,
 for any nontrivial Brauer  $(b,G)\-$pair $(Q,f)$ contained in $(P,e)\,,$ the block $f$ of $C_G (Q)$ admits $C_P (Q)$ as a defect group and it is inertially controlled.\/}
 \medskip
 \noindent
 {\bf Proof:} Firstly assume that $b$ is inertially controlled; let $(Q,f)$ be a Brauer
 $(b,G)\-$pair contained in  $(P,e)$ and choose a maximal Brauer  $\big(f,Q\.C_G (Q)\big)\-$ pair  $(R,g)\,;$ since $(Q,f)$ is also a
  Brauer  $\big(f,Q\.C_G (Q)\big)\-$pair, $(R,g)$ necessarily contains $(Q,f)$ (cf.~£1.6.1) and therefore
  it is also a Brauer  $(b,G)\-$pair; hence, there is $x\in G$ such that  (cf.~£1.6.1)
  $$(Q,f)^x\i (R,g)^x\i (P,e)
  \eqno £2.9.1\phantom{.}$$
  and therefore we get $x = zn$ for suitable $z\in C_G(Q)$ and $n\in N\,;$ so that the maximal Brauer
  $\big(f,Q\.C_G (Q)\big)\-$pair $(R,g)^z$ is contained in $(P,e)\,.$

  \smallskip
 Moreover, if $(T,h)$ is a  Brauer $\big(f,C_G (Q)\big)\-$pair, it is clear that
 $(Q\.T,h)$ is a   Brauer  $(b,G)\-$pair; conversely, by the argument above,
 $\big(C_P (Q),g^x\big)$ is a maximal Brauer $\big(f,C_G (Q)\big)\-$pair; then, if
 $\big(C_P (Q),g^x\big)$ contains $(T,h)$ and $(T,h)^z$
 with $z\in C_G(Q)\,,$ it is easily checked that $(P,e)$ contains $(Q\.T,h)$
 and~$(Q\.T,h)^z$ and therefore
 we still get $z = wn$ for suitable $w\in C_G(Q\.T)$ and~$n\in N\,,$ so that $n$ actually belongs to $C_N(Q)\,;$ consequently, since we have  $N/C_G(P)\cong L/C_L(P)\,,$ the block $f$ of $C_G (Q)$ is inertially controlled.

 \smallskip
 Conversely, arguing by contradiction, assume that $Q_\delta$ is an essential pointed group contained in~$P_\gamma\,.$ According to  [10,~Lemma~3.10], we may assume that the image of $N_P (Q)$ is a Sylow $p\-$subgroup
 of~$E_G (Q_\delta)$ and, since a proper subgroup $M$ of $E_G (Q_\delta)$ fulfilling condition~£2.8.1 above contains a Sylow $p\-$subgroup of  $E_G (Q_\delta)\,,$  we still may assume that $M$ contains the image
 of~$N_P(Q)\,.$ Moreover, it follows again from  [10,~Lemma~3.10] that there is
 a local pointed group $R_\varepsilon$ containing and normalizing $Q_\delta$ such that its image in~$E_G (Q_\delta)$ is not contained in $M\,;$ then, $R$ centralizes some nontrivial subgroup $Z$ of $Z(Q)$ and, denoting by $f$ the unique block of
  $H =C_G (Z)$ such that $(P,e)$ contains $(Z,f)\,,$ it follows from our hypothesis that $H\cap P$ is a defect group of this block.

 \smallskip
 Consequently, denoting by $h$ the block of $C_G (H\cap P)$ such that
$(P,e)$ contains $(H\cap P,h)\,,$ this pair is a maximal Brauer $(f,H)\-$pair;
moreover, $H$ contains $R$ and $C_G (Q)\,,$ and in particular we have
 $$(kH)(Q)\cong (kG)(Q)
 \eqno £2.9.2,$$
 so that ${\rm Br}_Q (\delta)$ determines a local point $\hat\delta$ of $Q$
 on $kH$ fulfilling
 $$E_H (Q_{\hat\delta}) \i E_G (Q_\delta)
 \eqno £2.9.3;$$
 then, appying again [10,~Lemma~3.10], we may assume that the image
 of $N_{H\cap P}(Q)$ in the intersection $E_H (Q_{\hat\delta})\cap M$  is a Sylow   $p\-$subgroup of $E_H (Q_{\hat\delta})\,,$  whereas this intersection does not contain the image of $R\,;$ hence, $Q_{\hat\delta}$ is an essential pointed group on $kHf\,,$ which contradicts our hypothesis. We are done.

 \bigskip
\bigskip
\noindent
{\bf £3. Strict semicovering homomorphism}
\bigskip

£3.1. Let $P$ be a finite $p\-$group, $B$ and $\hat B$ two $P\-$algebras and $g\,\colon B\to \hat B$
a unitary $P\-$algebra homomorphism; we say that $g$ is a {\it strict semicovering\/} if, for any
subgroup $Q$ of $P\,,$ we have ${\rm Ker}(g)^Q\i  J(B^Q)$ and the image $g(j)$ of a primitive idempotent
$j$ of $B^Q$ is still primitive in $\hat B^Q$ [6,~3.10]; namely if~$g$ induces a  homomorphism from the maximal semisimple quotient of~$B^Q$ to the  maximal semisimple quotient of~$\hat B^Q\,,$ mapping primitive idempotents
on primitive idempotents.
\eject

\medskip
£3.2. In other words, $g$ is a  strict semicovering
if and only if, for any subgroup $Q$ of $P\,,$ it induces a surjective map from the set of points of $Q$ on~$B$
to the set of points of $Q$ on $\hat B$ and, for any pair of mutually corresponding such points $\delta$ and
$\hat\delta\,,$ it induces a $k\-$algebra embedding [6,~3.10]
$$g(Q_\delta) : B(Q_{\delta})\too \hat B(Q_{\hat \delta})
\eqno £3.2.1.$$

\medskip
£3.3. Explicitly, if $g$ is a strict semicovering then, for any pointed group $Q_{\delta}$ on
$B\,,$ there is a unique point $\hat\delta$ of $Q$ on $\hat B$ fulfilling $g(\delta)\i \hat\delta\,;$ moreover,
this correspondence preserves {\it inclusion\/} and {\it localness\/} [6,~Proposition~3.15]. The composition of strict semicoverings is clearly a strict semicovering but, more precisely, the {\it strictness\/} provides a converse [6,~Proposition~3.6].

\bigskip
\noindent
{\bf Proposition~£3.4.} {\it  With the notation above, let $\hat g\,\colon \hat B\to \skew4\hat{\hat B}$
a second unitary $P\-$algebra homomorphism. Then, $\hat g\circ g$ is a strict semicovering if and only if
$\hat g$ and $g$ are so.\/}

\medskip
£3.5. The fact for a $P\-$algebra homomorphism of being a strict semico-vering is essentially of  ``local'' nature as it shows the following result [6,~Theorem~3.16].

\bigskip
\noindent
{\bf Theorem~£3.6.} {\it With the  notation above, the unitary $P\-$algebra homomorphism $g$ is a  strict semicovering if and only if, for any $p\-$subgroup $Q$ of~$P\,,$ the $\{1\}\-$algebra homomorphism
$$g(Q) : B(Q)\too \hat B(Q)
\eqno £3.6.1\phantom{.}$$
 induced by $g$ is a strict semicovering. \/}

 \medskip
 £3.7. Here, we may restrict ourselves to consider the following situation. Let $G$
 be a finite group, $H$ a normal subgroup of $G$ such that $G/H$ is a $p\-$group,
 $P$ a $p\-$subgroup of $G$ and $Z$ a subgroup of $Q = H\cap P$
 normal in $G$ and central in $H\,;$ set~$\bar G = G/Z$ and $\bar P = P/Z\,.$

\bigskip
\noindent
{\bf Proposition~£3.8.} {\it With the notation above, the canonical $\bar P\-$algebra homomorphism
$kH\to k\bar G$ is a semicovering.\/}
\medskip
\noindent
{\bf Proof:} For any subgroup $\bar Q = Q/Z$ of $\bar P\,,$ we have (cf.~£1.3.1)
$$(kH)(\bar Q)\cong kC_H (Q)\qq (k\bar G)(\bar Q) \cong kC_{\bar G}(\bar Q)
\eqno £3.8.1;$$
thus, a $p'\-$subgroup $K$ of the converse image  of $C_{\bar G}(\bar Q)$  centralizes $Q$ [5~Ch.~5, Theorem~3.2] and therefore it is contained in $C_H(Q)\,;$ that is to say, setting  $\overline{C_H(Q)}= C_H(Q)/Z\,,$
the quotient $C_{\bar G}(\bar Q)/\overline{C_H(Q)}$ is a $p\-$group.

\smallskip
Then, it follows from Lemma~£3.9 below that any simple $kC_{\bar G}(\bar Q)\-$mo-dule~$M$ has the form
$$M\cong {\rm Ind}_{kC_{\bar G}(\bar Q)_N}^{kC_{\bar G}(\bar Q)} (\hat N)
\eqno £3.8.2\phantom{.}$$
where $N$ is a simple $k\overline{C_H(Q)}\-$module, $kC_{\bar G}(\bar Q)_N$ the stabilizer
in $kC_{\bar G}(\bar Q)$ of the isomorphism class of $N$ and $\hat N$ the extended $kC_{\bar G}(\bar Q)_N\-$module. Moreover, any simple
$kC_H(Q)\-$module is also a simple $k\overline{C_H(Q)}\-$module and it appears in some simple
$kC_{\bar G}(\bar Q)\-$module. All this amounts to saying that the canonical $\{1\}\-$algebra homomorphism
$$kC_H (Q)\too kC_{\bar G}(\bar Q)
\eqno £3.8.3\phantom{.}$$
induces a homomorphism between the corresponding semisimple quotients preserving primitivity and then it suffices to apply Theorem~£3.6.

\bigskip
\noindent
{\bf Lemma~£3.9.} {\it  Let $X$ be a finite group and $Y$ a normal subgroup of $X$ such that $X/Y$ is a
$p\-$group. Then, any simple $kY\-$module $N$ can be extended to the stabilizer $X_N$ in $X$ of the isomorphism class of $N$ and, denoting by $\hat N$ the extended $kX_N\-$module,  ${\rm Ind}_{X_N}^X (\hat N)$ is a simple
$kX\-$module. Moreover, all the simple $kX\-$modules have this form.\/}
\medskip
\noindent
{\bf Proof:} Straightforward.

\bigskip
\noindent
{\bf Corollary~£3.10.} {\it With the same notation, let $\alpha = \{b\}$ be a point of $G$ on $kH$ and assume that
$P_\gamma$ is a defect pointed group of $G_\alpha\,;$ denote by $\bar b$ and $\bar\gamma$ the respective images in $k\bar G$ of~$b$ and $\gamma\,.$ Then, $b$ and $\bar b$ are respective blocks of $G$ and $\bar G\,,$ $\gamma$ and $\bar\gamma$ are respectively contained in local points $\tilde\gamma$ and $\tilde{\bar\gamma}$ of $P$
and~$\bar P$ on $kG$ and $k\bar G\,,$ and moreover $P_{\tilde\gamma}$  and~$\bar P_{\tilde{\bar\gamma}}$ are respective defect pointed groups of these blocks. In particular, setting $Q = H\cap P\,,$ $\bar H = H/Z$ and $\bar Q = Q/Z\,,$ the respective $P\-$ and $\bar P\-$interior algebras
$$(kH)_{\gamma}\otimes_Q P = \bigoplus_u (kH)_{\gamma}\. u\quad and\quad
(k\bar H)_{\bar\gamma}\otimes_{\bar Q}\bar P = \bigoplus_{\bar u} (k\bar H)_{\bar\gamma}\.\bar u
\eqno £3.10.1,$$
where $u\in P$  runs over a set of representatives for $P/Q$ and $\bar u$ is the image in $\bar P$ of $u\,,$
are respective source algebras of these blocks.\/}
\medskip
\noindent
{\bf Proof:} Since any block of $G$ is a $k\-$linear combination of $p'\-$elements of $G\,,$ $kH$ contains
all the blocks of $G$ and therefore $b$ is primitive in $Z(kG)\,;$ moreover, it is easily checked that $(kH)^G$
maps surjectively onto $(k\bar H)^{\bar G}$ and therefore $\bar\alpha = \{\bar b\}$ is also a point of $\bar G$
on $k\bar H\,,$ so that $\bar b$ is a block of $\bar G\,.$

\smallskip
Moreover, it follows from Propositions~£3.4 and~£3.8 that the canonical $\bar P\-$algebra homomorphisms
$$kH\too kG\qq kH\too  k\bar G
\eqno £3.10.2\phantom{.}$$
are strict semicovering; hence, $\gamma$ is contained in a local point $\tilde\gamma$ of $P$ on $kG$
and $\bar\gamma$ in a local point $\tilde{\bar\gamma}$ of $\bar P$ on $k\bar G\,;$ we claim that  $P_{\tilde\gamma}$  and~$\bar P_{\tilde{\bar\gamma}}$ are maximal local pointed groups on $kG$ and $k\bar G$ respectively.

\smallskip
Indeed, since the canonical homomorphism $kH\to kG$ is a semicovering,  a  local pointed group  $P'_{\tilde\gamma'}$
on $kG$ containing $P_{\tilde\gamma}$ comes from a local pointed group $P'_{\gamma'}$ on $kH$ and it is easily checked that $P'_{\gamma'}\i G_\alpha\,,$ so that we have $P'_{\gamma'}\i (P_\gamma)^x$ for a suitable $x\in G\,,$
which forces $P'_{\gamma'} = P_\gamma\,;$ since $\bar\alpha$ is a point of $\bar G$ on $k\bar H\,,$ the same argument proves that $\bar P_{\tilde{\bar\gamma}}$ is a maximal local pointed group on $k\bar G\,.$

\smallskip
The proof of the last statement is straightforward. We are done.

\bigskip
\bigskip
\noindent
{\bf £4. Stable embeddings: the proof of Theorem~£1.10}
\bigskip

£4.1. Let $G$ be a finite group and $A$ a $G\-$interior algebra; we say that a point $\beta$ of $H$ on $A$ is {\it projective\/}
if it is contained in $A^H_1$ or, equivalently, if it has a trivial defect group. Let $\hat A$ be a second  $G\-$interior algebra and
$f\,\colon \hat A\to A$ a $G\-$interior algebra homomorphism;  following [13,~6.4], we say that $f$ is a
{\it stable embedding\/}  if ${\rm Ker}(f)$ and $f(1_{\hat A}) A f(1_{\hat A})\big/f(\hat A)$ are projective  $k(G\times G)\-$modules or, equivalently, if the classe of the $k(G\times G)\-$module homomorphism
$$f : \hat A\too f(1_{\hat A}) A f(1_{\hat A})
\eqno £4.1.1\phantom{.}$$
 in the {\it stable category\/} of $k(G\times G)\-$modules is an isomorphism.

\medskip
£4.2. In this case, if $f$ is unitary, the exact sequence of $k(G\times G)\-$modules
$$0\too {\rm Ker}(f)\too \hat A\buildrel f\over\too A\too A/f(\hat A)\too 0
\eqno £4.2.1\phantom{.}$$
 is split [13,~6.4.1] and therefore, for any subgroup $H$ of $G\,,$ $f$ induces a
 {\it $C_G(H)\-$interior $N_G(H)\-$algebra isomorphism\/}
 $$\hat A^H/\hat A^H_1\cong A^H/A^H_1
 \eqno £4.2.2;$$
in particular, $f$ induces a bijection between the sets of {\it  nonprojective points\/} of $H$ on $\hat A$ and on $A$ and, for any pair of  corresponding nonprojective points $\hat\beta$~and $\beta\,,$ we have $N_G(H_{\hat\beta})
=N_G(H_\beta)\,,$ $f$ induces a $C_G(H)\-$interior
$N_G(H_\beta)\-$ algebra isomorphism [13,~4.6.2]
$$f(H_\beta) : \hat A(H_{\hat\beta})\cong A(H_{\beta})
\eqno £4.2.3\phantom{.}$$
and this isomorphism determines a central $k^*\-$extension isomorphism
$$\hat f(H_\beta) : \skew3\hat{\bar N}_G (H_{\hat\beta})\cong \skew3\hat{\bar N}_G (H_{\beta})
\eqno £4.2.4.$$
Moreover, this correspondence preserves {\it inclusion\/},  {\it localness\/} and
{\it fusions\/}.
\eject

 \medskip
 £4.3. We are ready to prove Theorem~£1.10; thus, $b$ is a block of $G\,,$   $P_\gamma$ is a defect pointed group of $b\,,$
  we set $N = N_G (P_\gamma)\,,$  $e$ is the corresponding block of $N\,,$ $\nu$ is the point of $N$ on $kG$ determined by $e\,,$ $\hat\gamma$ is  the local point of $P$ on $kN$ fulfilling ${\rm Br}_P(\hat\gamma) =  {\rm Br}_P(\gamma)$ and we denote by (cf.~£2.6.2)
 $$g : (kN)_{\hat\gamma}\too (kG)_{\gamma}
 \eqno £4.3.1\phantom{.}$$
 the unitary $P\-$interior algebra homomorphism determined as above by the multiplication by
 $\ell\in \nu\,;$ note that the restriction throughout $g$ induces a functor from the category of $kGb\-$modules to
the category of $kNe\-$modules which actually coincides with the functor determined by the $k(N\times G)\-$mo-dule~$\ell (kG)\,.$   Firstly, we prove a stronger form of the converse part.

\bigskip
\noindent
{\bf Proposition~£4.4.} {\it With the notation above, assume that the blocks $b$ of~$G$ and $e$ of $N$ are stably identical.
Then, for any nontrivial Brauer $(b,G)\-$pair $(Q,f)$ contained in $(P,e)\,,$ $N_P(Q)$ is a defect group of the block $f$
of $C_G(Q)\.N_P(Q)$ and a source algebra of  this block is isomorphic to $k_*\big(C_{\hat L}(Q)\.N_P(Q)\big)$ via an isomorphism inducing a $C_P(Q)\-$interior  algebra isomorphism from a source algebra of  the block $f$ of $C_G (Q)$  onto~$k_* C_{\hat L}(Q)\,.$\/}

\medskip
\noindent
{\bf Proof:} We can apply Theorem~6.9 and Corollary~7.4 in [13] to the Morita stable equivalences between $b$  and $e\,,$ and between $e$ and $b\,;$ in our present situation, $b$ and $e$ have the same defect group $P$ and, with the notation
in~[13], we may assume that $\ddot P = P$ and then $\ddot S = k$ is the trivial
$P\-$interior algebra and $\sigma = \sigma' = {\rm id}_P\,.$ Consequently, it follows from [13, 7.6.6] that the block $b$
of~$G$ is {\it inertially controled\/} and, for any nontrivial subgroup $Q$ of~$P\,,$ from~[13,~6.9.1] we get {\it $C_P(Q)\-$interior $N_P (Q)\-$algebra embeddings\/}
$$\eqalign{(kG)_\gamma (Q)\too (kN)_{\hat\gamma}(Q)\qq (kN)_{\hat\gamma} (Q)\too (kG)_{\gamma}(Q)\cr}
\eqno £4.4.1,$$
so that both are isomorphisms.

\smallskip
Now, we have  {\it $C_P(Q)\-$interior $N_P (Q)\-$algebra isomorphisms\/}
$$(kG)_{\gamma}(Q)\cong (kN)_{\hat\gamma} (Q)\cong k_*C_{\hat L}(Q)
\eqno £4.4.2\phantom{.}$$
and therefore the unity element is primitive in the $k\-$algebra
$$(kG)_{\gamma}(Q)^{C_P (Q)}\cong (kN)_{\hat\gamma}(Q)^{C_P (Q)}
\eqno £4.4.3;$$
thus, denoting by $f$ the block of $C_G(Q)$ such that $(Q,f)\i (P,e)$ [2,~Theorem~1.8], it is quite clear that
the $C_P(Q)\-$interior algebra $(kG)_{\gamma}(Q)$ is a source algebra of this block  and it is indeed isomorphic to
$k_*C_{\hat L}(Q)\,.$

\smallskip
Moreover,  it follows from Corollary~£3.10 above, applied to the groups $C_G(Q)\.N_P(Q)$ and $C_G(Q)\,,$ that
$N_P(Q)$ is a defect group of the block $f$ of $C_G(Q)\.N_P(Q)$ and that the $N_P (Q)\-$interior algebra
$$(kG)_\gamma (Q)\otimes_{C_P (Q)} N_P (Q) = \bigoplus_u\, (kG)_\gamma (Q)\.u
\eqno £4.4.4,$$
where $u\in N_P (Q)$ runs over a set of representatives for $N_P(Q)/C_P (Q)\,,$ is a source algebra
of this block; thus,  according to isomorphisms~£4.4.2, this  $N_P (Q)\-$interior algebra is isomorphic
to $k_*\big(C_{\hat L}(Q)\.N_P(Q)\big)\,.$
We are done.

 \bigskip
\noindent
{\bf Theorem~£4.5.} {\it With the notation above, for any nontrivial Brauer $(b,G)\-$pair $(Q,f)$ contained in $(P,e)\,,$ assume that $C_P(Q)$ is a defect group of the block~$f$ of $C_G (Q)$ and that a source algebra of this block is isomorphic to
$k_*\big(C_{\hat L}(Q)\big)\,.$ Then, $g$ is a stable embedding.\/}
\medskip
\noindent
{\bf Proof:}  Since $g$ is a {\it direct injection\/} of $k(P\times P)\-$modules
(cf.~£2.6), we have  ${\rm Ker}(g) =\{0\}$ and the quotient
$$M =(kG)_{\gamma}\big/g\big((kN)_{\hat\gamma}\big)
\eqno £4.5.1\phantom{.}$$
 is a direct summand of $(kG)_{\gamma}$ as $k(P\times P)\-$modules; hence,
 since  $(kG)_{\gamma}$ is a permutation $k(P\times P)\-$module, it suffices to
 prove that $M(W) =\{0\}$ for any nontrivial subgroup $W$ of $P\times P\,.$
 Actually, we have $(kG)(W) =\{0\}$ unless
 $$W = \Delta_\varphi (Q) = \{\big(u,\varphi (u)\big)\}_{u\in Q}
 \eqno £4.5.2\phantom{.}$$
 for some subgroup $Q$ of $P$ and some group homomorphism $\varphi\,\colon
 Q\to P$ induced by the conjugation by some $x\in G\,.$

 \smallskip
  More precisely, choosing  $i\in \gamma\,,$ the multiplication by $x$ on the right determines a $k\-$linear isomorphism
 $$(kG)_\gamma \big(\Delta_\varphi (Q)\big)\cong\big(i (kG)ix\big)(Q)
 \eqno £4.5.3;$$
thus, denoting by $f$ the block of $C_G (Q)$ such that $(P,e)$ contains
 $(Q,f)$ or, equivalently, such that $f{\rm Br}_Q(i) \not= 0\,,$ if we have $(kG)_\gamma \big(\Delta_\varphi (Q)\big)\not= \{0\}\,,$  we still have
  $f{\rm Br}_Q(i^x) \not= 0$  or, equivalently,  $(P,e)^x$ contains $(Q,f)$
  which amounts to saying that $\varphi\,\colon Q\to P$ is an
  $\F_{\!(b,G)}\-$morphism (cf.~£2.9). Hence, it suffices to prove that,
  for any nontrivial subgroup $Q$ of $P$ and any $\F_{\!(b,G)}\-$morphism
 $\varphi\,\colon Q\to P\,,$ we have $M \big(\Delta_\varphi (Q)\big) = \{0\}\,;$
 but, always  since $g$ is a {\it direct injection\/} of $k(P\times P)\-$modules,
 $g$ induces an injective homomorphism
$$g\big(\Delta_\varphi (Q)\big) : (kN)_{\hat\gamma}\big(\Delta_\varphi (Q)\big)\too (kG)_\gamma \big(\Delta_\varphi (Q)\big)
\eqno £4.5.4;$$
consequently, it suffices to prove that
$${\rm dim}\Big((kN)_{\hat\gamma}\big(\Delta_\varphi (Q)\big)\Big) = {\rm dim}\Big((kG)_\gamma \big(\Delta_\varphi (Q)\big)\Big)
\eqno £4.5.5\phantom{.}$$
and we argue by induction on $\vert P\colon Q\vert\,.$
\eject

\smallskip
Since we have a $P\-$interior algebra isomorphism $k_*\hat L\cong (kN)_{\hat\gamma}\,,$ we still have
$$(k_*\hat L)\big(\Delta_\varphi (Q)\big)\cong (kN)_{\hat\gamma}\big(\Delta_\varphi (Q)\big)
\eqno £4.5.6;$$
moreover, it is clear that $N_P(Q)$ centralizes a nontrivial subgroup $Z$ of $Z(Q)$
and then, according to our hypothesis, the $C_P(Z)\-$interior algebra
$k_* C_{\hat L}(Z) $ is isomorphic to a source algebra of the block $h$ of
$C_G (Z)$ such that $(P,e)$ contains $(Z,h)\,;$ in particular, setting $H = C_G(Z)\,,$
$(Q,f)$ is also a Brauer $(h,H)\-$pair, we have $C_H (Q) = C_G (Q)$ and $N_P (Q)$
remains a defect group of the block $f$ of $C_H (Q)\.N_P (Q)\,.$
Consequently, it easily follows from Proposition~£4.4 above, applied to the block $h$ of $H\,,$ that a source algebra
of the block $f$ of $C_G (Q)\.N_P (Q)$ is isomorphic to  $k_*\big(C_{\hat L} (Q)\.N_P(Q)\big)\,.$

\smallskip
At this point, we claim that in $(kG)_\gamma (Q)^{N_P(Q)}$ the unity element is primitive; since the point $\gamma$ is local, it follows from isomorphism~£2.2.2 that there is a primitive idempotent $\bar\ell$ of $(kG)_\gamma (Q)^{N_P (Q)}$ determining a local point of  $N_P (Q)$ on~$(kG)_\gamma (Q)\,;$ but, according to our induction hypothesis, for any subgroup~$R$ of~$N_P(Q)$ strictly containing $Q$ we may assume that that $(kN)_{\hat\gamma}(R)\cong
(kG)_\gamma (R)$ (cf.~£4.5.4) and, since ${\rm Br}_{\bar R}^{(kG)_\gamma (Q)}(\bar\ell) \not= 0$ where we set  $\bar R = R/Q$
(cf.~isomorphism~£2.2.2), we necessarily have
$${\rm Br}_{\bar R}^{(kG)_\gamma (Q)}(1_{(kG)_\gamma (Q)} -\bar\ell) = 0
\eqno £4.5.7;$$
thus,  the idempotent $1_{(kG)_\gamma (Q)} -\bar\ell$ belongs to [2, Lemmas~1.11 and~1.12]
 $$ \bigcap _R {\rm Ker}\big( {\rm Br}_{\bar R}^{(kG)_\gamma (Q)}\big)= \big((kG)_\gamma (Q)\big)_Q^{N_P(Q)} = {\rm Br}_Q\Big(\big((kG)_\gamma\big)_Q^P\Big)
\eqno £4.5.8\phantom{.}$$
 where $R$ runs over thet set of subgroups of~$N_P(Q)$ strictly containing $Q\,;$
but $0$ is the unique idempotent in $\big((kG)_\gamma\big)_Q^P\,;$
hence, we get $1_{(kG)_\gamma (Q)} =\bar\ell\,,$ proving the claim.

\smallskip
 Consequently, it follows from Corollary~£3.10 above, applied to the groups $C_G (Q)\.N_P (Q)$ and $C_G (Q)\,,$ that the
 $N_P (Q)\-$interior algebra (cf.~£2.3)
$$(kG)_\gamma (Q)\otimes_{C_P (Q)} N_P (Q) =
\bigoplus_u\, (kG)_\gamma (Q)\.u
\eqno £4.5.9,$$
where $u\in N_P (Q)$ runs over a set of representatives for $N_P(Q)/C_P (Q)\,,$ is a source algebra
of the block $f$ of $C_G (Q)\.N_P(Q)\,;$ hence, according to
 our hypothesis, we have a $N_P (Q)\-$interior algebra isomorphism
 $$k_*\big(C_{\hat L}(Q)\.N_P(Q¡\big)\cong (kG)_\gamma (Q)\otimes_{C_P (Q)} N_P (Q)
 \eqno £4.5.10;$$
now, according to isomorphism~£4.5.6 and equality~£4.5.9, we actually get
 $${\rm dim}\big((kN)_{\hat\gamma}(Q)\big) = {\rm dim}\big(k_*C_{\hat L}(Q)\big)
 = {\rm dim}\big((kG)_\gamma (Q)\big)
 \eqno £4.5.11\phantom{.}$$
 and therefore $g(Q)$ is an isomorphism.
 \eject

 \smallskip
  In particular, the interior $C_P (Q)\-$algebra $(kG)_\gamma (Q)
 \cong k_*C_{\hat L}(Q)$ is actually a source algebra of the block $f$ of $C_G (Q)$ and therefore, since
 we have (cf.~£1.11.1)
 $$C_{\hat L}(Q)\cong C_P(Q)\rtimes C_{\hat E}(Q)
 \eqno £4.5.12,$$
it follows from equalities~£2.7.2 that there is {\it no\/} essential pointed groups on~$kC_G (Q)f\,,$
so that the block $f$ of $C_G(Q)$ is {\it inertially controlled\/} (cf.~£2.9); hence,
it follows from Lemma~£2.9 and from our hypothesis that  the block~$b$ of~$G$ is also   {\it inertially controlled\/}.

\smallskip
Consequently, the $\F_{\!(b,G)}\-$morphism $\varphi\,\colon Q\to P$ above
is induced by some element $n\in N$ and therefore there is an inversible element
$a\in (kG)^P$ fulfil-ling $i^n = i^a\,,$ so that the multiplication by $na^{-1}$ on the right still determines a $k\-$linear isomorphism
 $$(kG)_\gamma \big(\Delta_\varphi (Q)\big)\cong (kG)_\gamma(Q)
 \eqno £4.5.13;$$
 similarly, we also get
 $$(kN)_{\hat\gamma} \big(\Delta_\varphi (Q)\big)\cong (kN)_{\hat\gamma}(Q)
 \eqno £4.5.14;$$
finally, equality~£4.5.5 follows from these isomorphisms and equality~£4.5.11.

\bigskip
\noindent
{\bf Corollary~£4.6.} {\it With the notation above, for any nontrivial Brauer $(b,G)\-$ pair $(Q,f)$ contained in $(P,e)\,,$ assume that $C_P(Q)$ is a defect group of the block~$f$ of $C_G (Q)$ and that a source algebra of this block is isomorphic to
$k_*\big(C_{\hat L}(Q)\big)\,.$ Then, the restriction throughout $g$ induces a stable equivalence between the categories of
$(kG)_\gamma\-$ and $(kN)_{\hat\gamma}\-$modules. In particular,
the blocks~$b$ of $G$ and $e$ of $N$ are stably identical.\/}
 \medskip
\noindent
{\bf Proof:}  With the notation in~£4.3 above, the indecomposable $k(N\times
G)$-module $\ell(k G)$ defined by the left-hand and the right-hand multiplication has
the $p\-$group $\Delta(P)=\{(u, u)|u\in P\}$ as a vertex and the trivial $k
\Delta(P)$-module $k$ as a source. Then this corollary follows from
Theorem~£4.5 above and [13,~Theorem~6.9] applied to the case where
$\ddot{M}=\ell(k G)$, $b=e$, $b'=b$, $P_\gamma=P_{\hat\gamma}$,
$P'_{\gamma'}=P_\gamma$ and $\ddot S=k\,.$

\bigskip
\bigskip
\noindent
{\bf £5. An inductive context: the proof of Theorem~£1.8}
\bigskip

£5.1. Let $G$ be a finite group, $b$ a block of $G$ and $P_\gamma$ a defect pointed group of $b\,;$
with the notation in £1.5 above, consider the following condition
\smallskip
\noindent
£5.1.1. {\it The block $b$ of $G$ is inertially controlled and, for any Brauer $(b,G)\-$pair $(Q,f)$ contained in $(P,e)\,,$ $f$ is a block of $C_G (Q)\.N_P (Q)$ with  trivial source simple modules.\/}
\smallskip
\noindent
First of all, we claim that if the block $b$ of $G$ fufills this condition then, for any Brauer $(b,G)\-$pair $(R,h)$ contained
in $(P,e)\,,$ the block $h$ of the group   $ H = C_G (R)$ fulfills the corresponding condition.
\eject

\medskip
£5.2. Indeed, it follows from Lemma~£2.9 that the block $h$ of $H$ is iner-tially controlled and that $T =C_P(R)$ is a defect group of  the block  $h$ of~$H\,;$ thus, denoting by $\ell$ the block  of~$C_G\big(R\.T\big)$ such that
$\big(R\.T,\ell\big)\i (P,e)$ [2,~Theorem~1.8], $(T,\ell)$ is a maximal Brauer $(h,H)\-$pair and,
if $(Q,f)$ is a Brauer  $(h,H)\-$pair contained in $(T,\ell)\,,$ $(R\.Q,f)$ is a Brauer $(b,G)\-$pair still contained in
$\big(R\.T,\ell\big)\i (P,e)$ and therefore $f$ is a block of  $C_G (R\.Q)\.N_P(R\.Q)$ with trivial source simple modules. Then, since  $C_H(Q)\.N_{T}(Q)$ is clearly subnormal~in~$C_G (R\.Q)\.N_P(R\.Q)\,,$
 it follows from  Lemma~£3.9, possibly applied  more than once, that $f$ is still a block of  $C_H(Q)\.N_{C_P(R)}(Q)$
with trivial source simple modules.

\medskip
£5.3.  At this point, assuming that the block $b$ of $G$ fufills condition~£5.1.1 and that, for any nontrivial Brauer $(b,G)\-$pair $(Q,f)$ contained in $(P,e)\,,$ we have  $C_G (Q)\not= G\,,$ it suffices to argue by induction
on $\vert G\vert$ to get the hypothesis of Theorem~£4.5, namely to get that, for any nontrivial Brauer $(b,G)\-$pair $(Q,f)$ contained in~$(P,e)\,,$ $C_P(Q)$ is a defect group of the block~$f$ of~$C_G (Q)$ (cf.~Lemma~£2.9) and that a source algebra of this block is isomorphic to~$k_*\big(C_{\hat L}(Q)\big)\,.$

\medskip
£5.4. In this situation, it follows from this theorem and from [13,~Theorem~6.9] that the blocks $b$ of~$G$ and $e$ of $N$ are  {\it stably identical\/} (cf.~£1.4);  more precisely, if $M$ is a simple $kGb\-$module of vertex $Q\i P$ and $f$ is the block of $C_G(Q)$ such that $(P,e)$  contains~$(Q,f)\,,$ on the one hand it follows from [8, Proposition~1.6] that the Brauer
$(b,G)\-$pair $(Q,f)$ is selfcentalizing, so that
$C_P (Q) = Z(Q)$ [16,~4.8 and~Corollary~7.3] and, on the other hand,  it easily follows from Theorem~£4.5 that the
$kNe\-$module  $\ell\.M\,,$ which is actually  indecomposable [7,~Theorem~2.1], has also vertex~$Q\,;$ moreover, since we are assuming that the trivial  $kQ\-$module $k$ is a source of~$M\,,$ it is clear that  the trivial $kQ\-$module $k$ is also a source
of~$\ell\. M\,.$

\medskip
£5.5.  Then, it follows again from [8,~Proposition~1.6] applied to the $N\-$interior algebra ${\rm End}_k(\ell\.M)\,,$
that the quotient $N_N(Q)/Q\,,$ and therefore the quotient $N_N(Q)/Q\.C_N (Q)$ [15,~Theorem~3.6], admit  blocks of {\it defect zero\/} --- namely, with trivial defect groups ---  which forces  [16,~1.19]
$$\Bbb O_p \big(N_N(Q)/Q\.C_N (Q)\big) =\{1\}
\eqno £5.5.1;$$
 but we have [11,~Proposition~14.6]
 $$C_P (Q) = Z(Q)\qq (kN)_{\hat\gamma}\cong k_*\hat L = k_*(P\rtimes \hat E^\circ)
 \eqno £5.5.2;$$
hence, denoting by $\hat\delta$ the unique local point of $Q$ on $kNe$ such that $P_{\hat \gamma}$ contains~$Q_{\hat \delta}$ (cf.~£2.8), it follows from~£2.7.2 and from the isomorphism in~£5.5.2 that, as in £1.11.1, we get [5,~Ch.~5, Theorem~3.4]
$$\eqalign{N_N(Q)/Q\.C_N (Q)&\cong E_N(Q_{\hat\delta})\cr
& = F_{(kN)_{\hat\gamma}}(Q_{\hat\delta})
\cong \big(N_P (Q)/Q\big)\rtimes N_{\hat E^\circ}(Q)\cr}
\eqno £5.5.3\phantom{.}$$
and, since $\Bbb O_p \big(N_N(Q)/Q\.C_N (Q)\big) =\{1\}\,,$ we still get $N_P (Q) = Q$ which forces $P = Q\,.$
In conclusion, $\ell\.M$ admits $P$ as a vertex and it has a trivial source, so that it is a simple $kNe\-$module according again to isomorphism~£5.5.2.

\medskip
£5.6. Finally, since the {\it stable equivalence\/} induced by the restriction throughout $g$ (cf.~£4.3.1) sends any simple $kGb\-$module to a simple $kNe\-$mo-dule, it follows from [7, Proposition~2.5] that  the restriction throughout $g$ actually
induces an equivalence of categories; moreover, since this equivalence is defined by a $k(G\times N)\-$module admitting a
$P\times P\-$stable basis (cf.~£4.3), it follows from [13,~Corollary~7.4 and~Remark~7.5] that the source algebras of the blocks $b$ of $G$ and $e$ of $N$ are isomorphic.

\medskip
£5.7. Assume now that the block $b$ of $G$ fufills condition~£5.1.1 and that there is an Abelian subgroup $Z$ of  $P$ such that $G = C_G (Z)\,;$ we are in the situation  considered in~£3.7 above with $H = G\,;$ hence, it follows from  Corollary~£3.10
that $\bar b$ is a block of~$\bar G$ and that $\bar\gamma$ is contained in a local point $\tilde{\bar\gamma}$ of $\bar P$
on $k\bar G$ such that $\bar P_{\tilde{\bar\gamma}}$ is a defect pointed group of $\bar b\,;$ denote by $\bar e$
the block of $C_{\bar G}(\bar P)$ determined by the point~$\tilde{\bar\gamma}\,.$

\medskip
£5.8. We claim that the block $\bar b$ of $\bar G$ fulfills the corresponding condition~£5.1.1. Indeed, if $(\bar Q,\bar f)$ is a Brauer $(\bar b,\bar G)\-$pair contained in $(\bar P,\bar e)$ and $Q$ is the converse  image of $\bar Q$ in $G\,,$ the image of $C_G(Q)$ in $C_{\bar G}(\bar Q)$ is a normal subgroup and, once again, the corresponding quotient is a $p\-$group  [5,~Ch.~5, Theorem~3.4]; hence, it follows again from Corollary~£3.10 that $\bar f$ is the image in $kC_{\bar G}(\bar Q)$
of a block  $f$ of the converse image $C$ of $C_{\bar G}(\bar Q)$ in~$G$
and then, since $C_G (Q)$ is normal in $C\,,$ it is quite clear that $f = {\rm Tr}_{C_{\tilde f}}^{C}(\tilde f)$ for a suitable block $\tilde f$ of~$C_G (Q)$ where $C_{\tilde f}$ denotes the stabilizer of $\tilde f$ in $C\,.$

\medskip
£5.9. More precisely, we claim that we can choose $\tilde f$ in such a way that $(P,e)$ contains $(Q,\tilde f)\,;$ indeed,
since $(\bar P,\bar e)$ contains $(\bar Q,\bar f)\,,$ there is a local point $\skew2\tilde{\bar\delta}$ of $\bar Q$ on $k\bar G$ such that we have $b_{\skew2\tilde{\bar\delta}} = \bar f$ and that $\bar P_{\tilde{\bar\gamma}}$ contains $Q_{\skew2\tilde{\bar\delta}}\,;$
then, it follows easily from Proposition~£3.8 and  from the obvious commutative diagram
$$\matrix{k\bar G^{\bar P}&\too&k\bar G^{\bar Q}\cr
\uparrow&&\uparrow\cr
k G^{ P}&\too&k G^{ Q}\cr}
\eqno £5.9.1\phantom{.}$$
that there is a point $\delta$ of $Q$ on $kG$ such that $P_\gamma$ contains $Q_\delta$ and that the image $\bar\delta$
of $\delta$ in $k\bar G$ is contained in $\skew2\tilde{\bar\delta}\,,$ which forces $\delta$ to be local;
at this point, it is easily checked that we can choose $\tilde f = b_\delta\,.$

\medskip
£5.10. Now, for any $\bar x\in \bar G$ such that $(\bar Q,\bar f)^{\bar x}\i (\bar P,\bar e)\,,$ the same argument
proves that we have $(Q,\tilde f)^{c x}\i (P,e)$ for some $x\in G$ lifting $\bar x$ and a suitable element $c$ of $C\,;$ then,
 since the block  $b$ of $G$ is inertially controlled, there are $n\in N$ and $z\in C_G(Q)$ fulfilling $cx = zn$ (cf.~£1.7)
 and therefore we get $\bar x = \bar c^{-1}\bar z\bar n$ where $\bar c\,,$ $\bar z$ and $\bar n$ denote the respective
  images of $c\,,$ $ z$ and $ n$ in~$\bar G\,,$ $ \bar c^{-1}\bar z$ centralizes $\bar Q$ and $\bar  n$ normalizes
  $(\bar P,\bar e)\,.$ This proves that the block $\bar b$ of~$\bar G$ is also inertially controlled.

\medskip
£5.11. Moreover, since $(Q,\tilde f)$ is a  Brauer $(b,G)\-$pair contained in $(P,e)\,,$ according to our hypothesis $\tilde  f$ is a block of  $C_G (Q)\.N_P(Q)$ with trivial source simple modules; but, since the block  $b$ of $G$ is inertially controlled, we have
$$E_G(Q,\tilde f)\cong  \big(N_P (Q)/Q\big)\rtimes N_{\hat E^\circ}(Q)
\eqno £5.11.1\phantom{.}$$
and therefore $C_{\tilde f}$ is contained in $C_G(Q)\.N_P(Q)\,;$ hence, since we have [2,~Theorem~1.8]
$$C_{\tilde f}\.N_P(Q) = C_G(Q)\.N_P(Q)\qq C_{\tilde f}\cap N_P(Q)  = C\cap N_P(Q)
\eqno £5.11.2,$$
 we clearly have [13,~2.6.4]
$$k\big(C\.N_P(Q)\big)  f \cong {\rm Ind}_{C_G (Q)\.N_P(Q)}^{C\.N_P(Q)}\Big(k\big(C_G(Q)\.N_P(Q)\big)\tilde f\Big)
\eqno £5.11.3\phantom{.}$$
and therefore  $f$ is also a block of  $C\.N_P(Q)$ with trivial source simple modules. Finally, since the $k\-$algebra
$k\big(C_{\bar G}(\bar Q)\.N_{\bar P}(\bar Q)\big)\bar f$ is the image of $k\big(C\.N_P(Q)\big) f\,,$
$\bar f$ is  a block of $C_{\bar G}(\bar Q)\.N_{\bar P}(\bar Q)$  with trivial source simple modules too.

\medskip
£5.12. Consequently, setting $\skew1\hat{\bar L} = \hat L/Z\,,$ it follows from our induction hypothesis that the source algebra of the block $\bar b$ of $\bar G$ is isomorphic to $k_*\skew1\hat{\bar L}$
and, in particular, we have
$${\rm dim}\big((k\bar G)_{\tilde{\bar\gamma}}\big) = \vert L\vert/\vert Z\vert
\eqno £5.12.1;$$
but, since the point $\tilde{\bar\gamma}$ contains the image of $\gamma\,,$ we may assume that
$(k\bar G)_{\tilde{\bar\gamma}}$ is the image of $(kG)_\gamma$ or, equivalently, that
$$(k\bar G)_{\tilde{\bar\gamma}}\cong k\otimes_{kZ} (kG)_\gamma
\eqno £5.12.2\phantom{.}$$
and, in particular, we get
$${\rm dim}\big((kG)_{\gamma}\big) = \vert Z\vert\, {\rm dim}\big((k\bar G)_{\tilde{\bar\gamma}}\big) =
 \vert L\vert
 \eqno £5.12.3;$$
 hence, the unitary $P\-$interior algebra homomorphism~£2.6.2 is actually an isomorphism
 $$k_*\hat L\cong (kN)_{\hat\gamma}\cong (kG)_{\gamma}
\eqno £5.12.4.$$

\medskip
£5.13. Conversely, assume that the source algebra $(kG)_\gamma$ is isomorphic to~$k_*\hat L\,,$
so that the unitary $P\-$interior algebra homomorphism~£2.6.2 is  an isomorphism; then, it follows from
equalities~£2.7.2 applied to the blocks $b$ of~$G$ and $\{1\}$ of $\hat L$ that there are no essential
pointed groups on $kGb$ (cf.~£2.8) and therefore the block $b$ of $G$ is inertially controled (cf.~£2.9).
 \eject

 \medskip
 £5.14. For any Brauer $(b,G)\-$pair $(Q,f)$ contained in $(P,e)\,,$ since we have (cf.~£1.3.1 and~£1.11.1)
 $$\eqalign{(kG)(Q)\cong kC_G (Q)\qq
 (kG)_\gamma(Q) &\cong (k_*\hat L)(Q) \cong k_*C_{\hat L}(Q)\cr
 & = k_*\big(C_P (Q)\rtimes C_{\hat E}(Q)\big)\cr}
 \eqno £5.14.1,$$
 the $C_P (Q)\-$interior algebra $(kG)_\gamma(Q)$ is a source algebra of the block $f$ of~$C_G (Q)\,;$
 then, it follows from Corollary~£3.10 that a source algebra of the block $f$ of
 $C_G (Q)\.N_P (Q)$ is isomorphic to the $N_P (Q)\-$interior algebra
 $$(kG)_{\gamma}(Q)\otimes_{C_P (Q)} N_P(Q)
 \eqno £5.14.2;$$
finally, according to isomorphisms~£5.14.1, this  $N_P (Q)\-$interior algebra is isomorphic to
 $$(k_*\hat L)(Q)\otimes_{C_P (Q)} N_P(Q)\cong k_*N_{\hat L}(Q)=   k_*\big(N_P (Q)\rtimes N_{\hat E}(Q)\big)
 \eqno £5.14.3\phantom{.}$$
 which clearly has trivial source simple modules. We are done.

\bigskip
\bigskip

\noindent{\bf Acknowledgement}\bigskip

When preparing the Abelian defect group case of this paper, the
second author was supported by the Alexander von Humboldt Foundation
of Germany and stayed at Jena University with the host Professor
Burkhard K\"ulshammer. He thanks the Alexander von Humboldt
Foundation very much for its support and Professor Burkhard
K\"ulshammer and his family for their hospitality. The second author
is also supported by the Key Project of Chinese Ministry of
Education and Program for New Century Excellent Talents in
University and he thanks the two supports a lot.

\bigskip
\bigskip

\noindent
{\bf References}
\bigskip
\noindent
[1]\phantom{.} Jon Alperin,  {\it Weight for finite groups\/}, in Proc. Symp.
Pure Math. 47(1987) 369-379, Amer. Math. Soc., Providence.
\smallskip\noindent
[2]\phantom{.} Michel Brou\'e and Lluis Puig, {\it Characters and Local
Structure in $G\-$al-gebras,\/} Journal of Algebra, 63(1980), 306-317.
\smallskip\noindent
[3]\phantom{.} Susanne Danz and Burkhard K\"ulshammer, {\it Vertices, sources and Green correspondents
of the simple modules for the large Mathieu groups,\/} Journal of Algebra, 322(2009) 3919-3949.
\smallskip\noindent
[4]\phantom{.} Yun Fan and Lluis Puig, {\it On blocks with nilpotent
coefficient extensions\/}, Algebras and Representation Theory, 1(1998),
27-73 and Publisher revised form, 2(1999), 209.
\smallskip\noindent
[5]\phantom{.} Daniel Gorenstein, {\it ``Finite groups''\/} Harper's Series,
1968, Harper and Row.
\smallskip\noindent
[6]\phantom{.} Burkhard K\"ulshammer and Llu\'\i s Puig, {\it Extensions of
nilpotent blocks}, Inventiones math., 102(1990), 17-71.

\smallskip\noindent
[7]\phantom{.} Markus Linckelmann, {\it Stable equivalences of Morita type for self-injective
algebras and $p\-$groups\/}, Math. Z. 223(1996), 87-100.
 \smallskip\noindent
[8]\phantom{.} Llu\'\i s Puig, {\it Pointed groups and  construction of
characters}, Math. Zeit. 176(1981), 265-292.
\smallskip\noindent
[9]\phantom{.} Llu\'\i s Puig, {\it Local fusions in block source algebras\/},
Journal of Algebra, 104(1986), 358-369.
\smallskip\noindent
[10]\phantom{.} Llu\'\i s Puig, {\it Nilpotent blocks and their source
algebras}, Inventiones math., 93(1988), 77-116.
\smallskip\noindent
[11]\phantom{.} Llu\'\i s Puig, {\it Pointed groups and  construction of
modules}, Journal of Algebra, 116(1988), 7-129.
\smallskip\noindent
[12]\phantom{.} Llu\'\i s Puig, {\it Alg\`ebres de source de certains blocks des groupes
de Chevalley\/}, in {\it ``Repr\'esentations lin\'eaires des groupes finis''\/}, Ast\'erisque, 181-182 (1990),
Soc. Math. de France
\smallskip\noindent
[13]\phantom{.} Llu\'\i s Puig, {\it ``On the Morita and Rickard
equivalences between Brauer blocks''\/}, Progress in Math., 178(1999), Birkh\"auser,
Basel.
\smallskip\noindent
[14]\phantom{.} Llu\'\i s Puig, {\it Source algebras of $p\-$central group
extensions\/}, Journal of Algebra, 235(2001), 359-398.
\smallskip\noindent
[15]\phantom{.} Llu\'\i s Puig, {\it Block Source Algebras  in p-Solvable
Groups},  Michigan Math. J. 58(2009), 323-328
\smallskip\noindent
[16]\phantom{.} Llu\'\i s Puig, {\it ``Frobenius categories versus Brauer blocks''\/}, Progress in Math.,
274(2009), Birkh\"auser, Basel.
\smallskip\noindent
[17]\phantom{.} Leonard Scott, {\it Defect groups and the isomorphism problem\/}, in {\it ``Repr\'e-sentations lin\'eaires des groupes finis''\/}, Ast\'erisque, 181-182 (1990),
Soc.Math. de France
\smallskip
\noindent
[18]\phantom{.} Yuanyang Zhou, {\it Gluing Morita equivalences induced by $p\-$permutation
modules\/}, preprint

\bigskip
\bigskip
\noindent
{\bf Abstract.} Motivated by an observation in [3], we determine the {\it source algebra\/}, and therefore all the structure, of the blocks without {\it essential Brauer pairs\/} where the simple modules of all the {\it Brauer corespondents\/} have trivial sources.

\end